\documentclass[12pt]{article}

\scrollmode

\usepackage{amsmath,amssymb,amsfonts,graphicx}
\usepackage[latin1]{inputenc}

\newtheorem{thm}{Theorem}[section]
\newtheorem{prop}[thm]{Proposition}

\newtheorem{defn}[thm]{Definition}

\newtheorem{rem}[thm]{Remark}

\def\be#1 {\begin{equation} \label{#1}}
\newcommand{\ee}{\end{equation}}

\def\sqw{\hbox{\rlap{\leavevmode\raise.3ex\hbox{$\sqcap$}}$%
\sqcup$}}
\def\findem{\ifmmode\sqw\else{\ifhmode\unskip\fi\nobreak\hfil
\penalty50\hskip1em\null\nobreak\hfil\sqw
\parfillskip=0pt\finalhyphendemerits=0\endgraf}\fi}

\newcommand{\mb}{\medskip\noindent}
\newcommand{\gb}{\bigskip\noindent}
\newcommand{\R}{\mathbb R}

\newcommand{\N}{\mathbb N}

\newcommand{\s}{\mathbf S}
\newcommand{\F}{\mathcal F}

\def\deme {\noindent {\bf Proof : }}

\begin{document}

\title{On Maximal $L^p$-regularity}

\author{ F. Bernicot\\ Universit\'e de Paris-Sud \\F-91405 Orsay
Cedex\\frederic.bernicot@math.u-psud.fr \and J. Zhao \\ Beijing Normal University  \\ Beijing 100875, P.R. China \\ jzhao@bnu.edu.cn}

\date{ }

\maketitle

\begin{abstract}
  The aim of this paper is to propose weak assumptions to prove maximal $L^q$ regularity for Cauchy problem~:
$$\frac{d u}{dt}(t)-Lu(t)=f(t).$$ 
Mainly we only require ``off-diagonal'' estimates on the real semigroup $(e^{tL})_{t>0}$ to obtain maximal $L^q$ regularity. The main idea is to use a one kind of Hardy space $H^1$ adapted to this problem and then use interpolation results. These techniques permit us to prove weighted maximal regularity too.
\end{abstract}

\mb {\bf Key-words~:}  Maximal regularity, Hardy spaces, atomic decomposition. 

\mb {\bf  AMS2000 Classification :} 42B25, 42B30, 42B35, 46M35
\
\tableofcontents

\section{Introduction}

\mb Let $(Y,d_Y,\nu)$ be a space of homogeneous type. Let $L$ be the
infinitesimal generator of an analytic semigroup of operators on
$L^p:=L^p(Y)$ and $J=(0,l],$ $l>0$ or $J=(0,+\infty)$ (in the second
case, one has to assume that $L$ generates a bounded analytic
semigroup).

\gb Consider the Cauchy problem
\begin{equation} \left\{\begin{array}{ll}
\frac{du}{dt}(t)-Lu(t)=f(t), &  t\in J,\\
u(0)=0, &
\end{array}\right. \tag{CP}
\end{equation}
where $f: J\rightarrow B$ is given, where $B$ is a Banach space. If
$e^{tL}$ is the semigroup generated by $L$, $u$ is formally given by
$$u(t)=\int_0^t e^{(t-s)L}f(s)ds.$$
For fixed $q \in (1,+\infty)$, one says that there is maximal $L^q$
regularity on $B = L^p$ for the problem if for every $f \in
L^q(J, L^p)$, $\frac{\partial u}{\partial t}$(or $Lu$) belongs to
$L^q(J, L^p)$. 
It is known that the property of maximal $L^q$-regularity does not depend on $q\in(1,\infty)$.

\mb For the maximal $L^q$ regularity, we refer the reader to the works of P. Cannarsa and V. Vespri \cite{CV}, T. Coulhon and X.T. Duong \cite{CD,CD1}, L. de Simon \cite{DS}, M. Hieber and J. Prss \cite{HP} and D. Lamberton \cite{L} etc. The literature is so vast that we do not give exhaustive
references. However we emphasize that in all these works, the different authors obtain maximal regularity under the assumption that the heat kernel (the kernel of the semigroup) admits pointwise estimates and gaussian decays. Such assumptions imply that the semigroup extends consistently to all Lebesgue spaces $L^p$ for $p\in(1,\infty)$. \\
For a few years, people have studied problems associated to a semigroup, which do not satisfy this property. For example, gaussian estimates have been succesfully generalized by ``off-diagonal estimates'' for studying the boundedness of Riesz transforms on a manifold (see \cite{A}). That is why we look for weaker assumptions associated to ``off-diagonal'' estimates on the semigroup to guarantee the maximal $L^q$-regularity.

\gb In this direction, there is a first work of S. Blunck and P.C. Kunstmann \cite{BK2}. The authors have obtained the following result (using the $R$-boundedness of the complex semigroup and the recent characterization of L. Weiss \cite{weis})~:

\begin{thm} \label{th} Let $\delta$ be the homogeneous dimension of $Y$. Assume that $(e^{zL})_z$ is a bounded analytic semigroup on $L^2$ and $p_0<2<q$ be  exponents. Suppose there are coefficients $(g(k))_{k\geq 1}$ such that for all balls $Q$ or radius $r_Q$ and all integer $k\geq 0$, we have
\be{assum1} \left\|{\bf 1}_{Q} e^{r_Q^2L} {\bf 1}_{(k+1)Q\setminus kQ} \right\|_{L^{p_0} \to L^q} \lesssim \nu(Q)^{\frac{1}{q} - \frac{1}{p_0} } g(k) \ee
with
\be{assum2} \sum_{k=1}^\infty k^{\delta-1} g(k) <\infty.\ee
Then for all $r\in(p_0,2]$, $L$ has maximal regularity on $L^r$.
\end{thm}

\gb Now we come to our results. We look for similar results with some improvements. First as the conclusion concerns only exponents $r\in(p_0,2]$, we would like to not require assumption (\ref{assum1}) with an exponent $q>2$. In addition, we want to understand how the assumption (\ref{assum2}) is important.
In our result, we will give some similar assumption, which seem to be not comparable to this one. However our proof (which use very different techniques) permit us to obtain simultaneously positive and new results for ``weighted'' maximal regularity.

\mb In \cite{ABZ}, the authors consider the Cauchy problem $(CP)$ with $-L$ equals to the Laplacian operator on some Riemannian manifolds or a sublapacian on some Lie groups or some second order elliptic operators on a domain.
We show the boundedness of the operator of maximal regularity
$f\mapsto Lu$ and its adjoint on appropriate Hardy spaces. In this
paper, we apply the general theories of our paper \cite{bz} to
the maximal regularity in abstract setting. In
\cite{bz}, we  construct
 Hardy spaces through an atomic (or molecular) decomposition which
  keep the main properties of the (already known) Hardy spaces $H^1$. We prove some results about continuity from these spaces into $L^1$ and some results about interpolation between these spaces and the Lebesgue spaces. Now we
  will use these theories to study the maximal regularity. 

\mb Here is our main result~:

\begin{thm} Let $L$ be a generator of a bounded analytic semigroup ${\mathcal T}:=(e^{tL})_{t>0}$ on $L^2(Y)$ such that
${\mathcal T},(tLe^{tL})_{t>0}$ and $(t^2L^2e^{tL})_{t>0}$ satisfy ``$L^2-L^2$ off-diagonal decay'' (precisely belonging to the class ${\mathcal O}_{4}(L^2-L^2)$, see Definition \ref{def:off}). For an exponent $p_0\in (1,2]$, we assume some weak ``$L^{p_0}-L^2$ off-diagonal decay'' (we require (\ref{hypi}), see Proposition \ref{assumcas}). \\
Then for all exponent $p\in(p_0,2]$, the operator $T$ admits a continuous extension on $L^p(Y)$ and so $L$ has maximal $L^p$-regularity.
In addition we can have weighted results~: let $\omega\in {\mathbb A}_\infty(Y)$ be a weight on $Y$. Then for all exponents $p\in(p_0,2)$ satisfying
$$ \omega\in {\mathbb A}_{p/p_0} \cap RH_{(2/p)'},$$
$L$ has maximal $L^p(\omega)$-regularity.
\end{thm}

\begin{rem} Our weak ``$L^{p_0}-L^2$ off-diagonal decay'' is similar to the assumption (\ref{assum1}) of Theorem \ref{th} but is not comparable to this one. What is important is that we only require informations on the real semigroup. In addition, the answer concerning weighted results is totally new and do not seem accessible by the techniques of \cite{BK2} used to prove Theorem \ref{th}.
\end{rem}

\mb The plan is as follows ~:  in Section \ref{section2}, we recall the abstract results concerning Hardy spaces. Then in Section \ref{section3}, we will explain the application to the maximal regularity problem : how to define an adapted Hardy space. Then we conclude in Section \ref{section4} by checking the abstract assumptions for this application. We will finish in the last section to give results for exponents $p\geq 2$ and study the Hardy spaces adapted to this problem of maximal regularity.

\section{Preliminaries} \label{section2}

In this section, we give an overview of some basic facts which we
will use in the sequel. For more details concerning abstract Hardy
spaces, see \cite{bz}.

\mb Let $(X,d,\mu)$ be a space of homogeneous type. That is meaning $d$ is a
quasi-distance on the space $X$ and $\mu$ a Borel measure which
satisfies the doubling property~: \be{homogene} \exists A>0, \
\exists \delta>0, \qquad \forall x\in X,\forall r>0, \forall t\geq
1,\qquad \mu(B(x,tr)) \leq A t^{\delta}\mu(B(x,r)), \ee where
$B(x,r)$ is the open ball with center $x\in X$ and radius $r>0$. We
call $\delta$ the homogeneous dimension of $X$.

\mb Let $Q$ be a ball, for $i\geq 0$, we write $S_i(Q)$ the scaled
corona around the ball $Q$~:
$$ S_i(Q):=\left\{ x,\ 2^{i} \leq 1+\frac{d(x,c(Q))}{r_Q} < 2^{i+1} \right\},$$
where $r_Q$  is the radius of $Q$ and $c(Q)$ is its center.
 Then  $S_0(Q)$ corresponds to the ball $Q$ and $S_i(Q) \subset 2^{i+1}Q$ for $i\geq 1$, where $\lambda Q$ is as usual the ball with center $c(Q)$ and radius $\lambda r_Q$.

\mb Denote $\mathcal{Q}$ the collection of all balls~: $ \mathcal{Q}:=
\left\{ B(x,r),\ x\in X, r>0 \right\},$ and $\mathbb{B}:=(B_Q)_{Q\in
\mathcal{Q}}$ a collection of $L^2$-bounded linear operators, indexed
by the collection $\mathcal{Q}$. We assume that these operators
$B_Q$ are uniformly bounded on $L^2$ : there exists a constant
$0<A'<\infty$ so that~: \be{operh} \forall f\in L^2 ,\ \forall  Q
\textrm{ ball}, \qquad \|B_Q(f) \|_2 \leq A'(\delta) \|f\|_2. \ee

\mb Now, we recall some definitions and theorems of \cite{bz}. The $\epsilon$-molecules (or atoms) are defined as follows.

\begin{defn}(\cite{bz}) Let $\epsilon>0$ be a fixed parameter.
A function $m\in L^{1}_{loc}$ is called an $\epsilon$-molecule associated to a ball $Q$ if there exists a real function $f_Q$ such that
$$m=B_Q(f_Q),$$
with
$$\forall i\geq 0, \qquad  \|f_Q\|_{2,S_i(Q)} \leq \left(\mu(2^{i}Q)\right)^{-1/2} 2^{-\epsilon i}.$$
We call $m=B_Q(f_Q)$ an atom if in addition we have $supp(f_Q) \subset Q$.
So an atom is exactly an $\infty$-molecule.
\end{defn}

\mb Using this definition, we can define the ``finite'' molecular (atomic) Hardy space. 

\begin{defn}(\cite{bz}) A measurable function $h$ belongs to the ``finite'' molecular Hardy space $H^1_{F,\epsilon,mol}$ if there exists a finite decomposition~:
$$h=\sum_{i} \lambda_i m_i  \qquad \mu-a.e, $$
where for all $i$, $m_i$ is an $\epsilon$-molecule and $\lambda_{i}$
are real numbers satisfying
$$\sum_{i\in \N} |\lambda_i| <\infty. $$
We define the norm~:
$$\|h\|_{H^1_{F,\epsilon,mol}}:= \inf_{h=\sum_{i} \lambda_i m_i} \sum_{i} |\lambda_i|,$$
where we take the infimum over all the finite atomic decompositions.
Similarly we define the ``finite'' atomic space $H^1_{F,ato}$ replacing $\epsilon$-molecules by atoms.
\end{defn}

\mb We will use the following theorem for studying maximal regularity.

\begin{prop} (\cite{bz}) \label{theo2h} Let $T$ be an $L^2$-bounded sublinear operator satisfying the following ``off-diagonal'' estimates~:
for all ball $Q$, for all $k\geq 0,j\geq 2$, there exists some coefficient $\alpha_ {j,k}(Q)$ such that for every $L^2$-function $f$ supported in $S_k(Q)$
 \be{assu} \left(\frac{1}{\mu(2^{j+k+1}Q)} \int_{S_j(2^kQ)} \left| T(B_Q(f))\right|^2 d\mu \right)^{1/2} \leq \alpha_{j,k}(Q) \left(\frac{1}{\mu(2^{k+1}Q)}\int_{S_k(Q)} |f|^2 d\mu \right)^{1/2}. \ee
If the coefficients $\alpha_{j,k}$ satisfy \be{hyph} \Lambda:= \sup_{k\geq 0}\  \sup_{Q \textrm{ ball}} \
\left[\sum_{j\geq 2} \frac{\mu(2^{j+k+1}Q)}{\mu(2^{k+1}Q)} \alpha_{j,k}(Q) \right] <\infty, \ee then for all $\epsilon>0$ there exists a constant $C=C(\epsilon)$ such that
$$ \forall f\in H^1_{F,\epsilon, mol} \qquad \|T(f)\|_{1} \leq C \|f\|_{H^1_{F,\epsilon, mol}}.$$
\end{prop}

\begin{defn} (\cite{bz}) We set $A_Q$ be the operator $Id-B_Q$. For $\sigma\in [2,\infty]$ we define the maximal operator~:
\be{opeM} \forall x\in X, \qquad M_{\sigma}(f)(x):= \sup_{\genfrac{}{}{0pt}{}{Q \textrm{ball}}{x\in Q}}\  \left( \frac{1}{\mu(Q)} \int_Q \left|A_Q^*(f)\right|^\sigma d\mu  \right) ^{1/\sigma}. \ee
We use duality so we write $A_Q^*$ for the adjoint operator.
The standard maximal ``Hardy-Littlewood'' operator is defined by ~: for $s>0$,
$$\forall x\in X, \qquad M_{HL,s}(f)(x):=\sup_{\genfrac{}{}{0pt}{}{Q \textrm{ball}}{x\in Q}} \left(\frac{1}{\mu(Q)} \int_{Q} \left|f \right|^{s} d\mu
\right)^{1/s}.$$
\end{defn}

\mb The main result of \cite{bz} is the following one about interpolation between $L^2$ and the Hardy spaces~:

\begin{thm} (\cite{bz}) \label{theogeneh}  Let $\sigma\in(2,\infty]$. Assume that  we have an implicit constant such that for all functions $h\in L^2$
$$ M_\sigma(h) \lesssim M_{HL,2}(h).$$
Let $T$ be an $L^2$-bounded, linear operator. Assume that $T$ is continuous from $H^1_{F,ato}$ (or $H^1_{F, \epsilon, mol}$) into $L^1$. Then for all exponent $p\in (\sigma',2]$ there exists a constant $C=C(p)$ such that~:
$$ \forall f\in L^2 \cap L^p,\qquad \|T(f)\|_{p}\leq C\|f\|_{p}.$$
\end{thm}

\mb We have boundedness in weighted spaces too.
We recall the definition of Muckenhoupt's weights and Reverse H\"older classes~:

\begin{defn} A nonnegative function $\omega$ on $X$ belongs to the class ${\mathbb A}_p$ for $1< p <\infty$ if
$$\sup_{Q \textrm{ ball}} \left(\frac{1}{\mu(Q)}\int_Q \omega d\mu \right) \left( \frac{1}{\mu(Q)}\int_Q \omega^{-1/(p-1)} d\mu \right)^{p-1} <\infty.$$
A nonnegative function $\omega$ on $X$ belongs to the class $RH_q$ for $1<q<\infty$, if there is a constant $C$ such that for every ball $Q\subset X$
$$ \left( \frac{1}{\mu(Q)}\int_Q \omega^q d\mu \right)^{1/q} \leq C \left(\frac{1}{\mu(Q)}\int_Q \omega d\mu \right).$$
We use the following notation of \cite{AMR}~: \\
Let $\omega\in {\mathbb A}_\infty$ be a weight on $X$ and $0<p_0<q_0\leq \infty$ be two exponents, we introduce the set
$$ \mathcal{W}_{\omega}(p_0,q_0):=\left\{ p\in(p_0,q_0),\ \omega\in {\mathbb A}_{p/p_0} \cap RH_{(q_0/p)'} \right\}.$$
\end{defn}

\mb Then we have the following result~:

\begin{thm} \label{theopoids} Let $\sigma\in(2,\infty]$. Assume that  we have an implicit constant such that for all $h\in L^2$
$$ M_\sigma(h) \lesssim M_{HL,2}(h).$$
Let $T$ be an $L^2$-bounded, linear operator such that for all balls $Q$ and for all functions $f$ supported in $Q$
\be{hypp} \forall j\geq 2 \quad \left(\frac{1}{\mu(2^{j+1}Q)} \int_{S_j(Q)} \left| T(B_Q(f))\right|^2 d\mu \right)^{1/2} \leq \alpha_j(Q) \left(\frac{1}{\mu(Q)}\int_{Q} |f|^2  d\mu \right)^{1/2}, \ee
with coefficients $\alpha_j(Q)$ satisfying
$$\sup_{Q \textrm{ ball}} \sum_{j\geq 0} \frac{\mu(2^{j+1}Q)}{\mu(Q)} \alpha_j(Q) <\infty.$$
Let $\omega\in {\mathbb A}_\infty$ be a weight. Then for all exponents $p\in {\mathcal W}_{\omega}(\sigma',2)$, there exists a constant $C$ such that
$$ \forall f\in L^2 \cap L^p(\omega), \qquad  \|T(f)\|_{p,\omega d\mu} \leq C \|f\|_{p,\omega d\mu}.$$
\end{thm}

\begin{rem} From Proposition \ref{theo2h}, (\ref{hypp}) implies the $H^1_{F,ato}-L^1$ boundedness of $T$. However, the proof for weighted results requires (\ref{hypp}) and not only this boundedness of the operator $T$.
\end{rem}

\mb Now we give some results concerning the Hardy spaces. Assume that ${\mathbb B}$ satisfies some decay estimates : for $M>n/2$ an integer (with $n$ the homogeneous dimension of $X$), there exists a constant $C$ such that
\be{decay} \forall i\geq 0,\ \forall k\geq 0,\ \forall f\in L^2,\
\textrm{supp}(f) \subset 2^k Q \qquad \left\| B_Q(f)
\right\|_{2,S_i(2^{k}Q)} \leq C 2^{-Mi} \|f\|_{2,2^{k}Q}. \ee

\mb Then we have the following results:

\begin{prop} (\cite{bz}) \label{contL1} The spaces $H^1_{ato}$ and $H^1_{\epsilon, mol}$ are Banach spaces.
And
$$ \forall \epsilon>0, \qquad H^{1}_{ato} \hookrightarrow H^{1}_{\epsilon, mol} \hookrightarrow L^1.$$
Therefore
$$L^\infty \subset (H^1_{\epsilon ,mol})^{*} \subset (H^1_{ato})^{*}.$$

\end{prop}

\mb We denote $H^1_{CW}$ the classical Hardy space (of Coifman-Weiss) (see \cite{CW}). As we noted in \cite{bz}, it corresponds to our Hardy space $H^1_{ato}$ or $H^1_{\epsilon, mol}$ when the operators $B_Q$ exactly correspond to the oscillation operators.

\begin{prop} (\cite{bz}) \label{inclus} Let $\epsilon\in(0,\infty]$. The inclusion $H^1_{\epsilon, mol} \subset H^1_{CW}$ is equivalent to the fact that for all $Q\in {\mathcal Q}$, $(A_Q)^* ({\bf 1}_X)={\bf 1}_{X}$ in $(Mol_{\epsilon,Q})^{*}$.
In this case for all $\epsilon'\geq \epsilon$ we have the inclusions $H^1_{ato} \subset H^1_{\epsilon',mol} \subset H^1_{\epsilon, mol} \subset H^1_{CW}$.
\end{prop}

\section{An application to maximal $L^q$ regularity on Lebesgue spaces.}
\label{section3}

In this section, we apply the previous general theory to maximal
$L^q$ regularity for Cauchy Problem.

\mb We first define an operator $T$~:
\begin{defn} With $L$ the generator of the semigroup, we define the operator~:
$$ Tf(t,x)=\int_0^t  \left[Le^{(t-s)L} f(s,.) \right](x) ds. $$
\end{defn}

\mb Let $p,q\in(1,\infty)$ be two exponents. We know that the maximal $L^q$ regularity on $L^p(Y)$ is equivalent to the fact that $T$ is bounded on $L^p(J \times Y)$. That is why we study this operator.
Of course, the problem of maximal $L^q$ regularity is completely understood by the abstract result in \cite{weis} of L. Weis using the $R$ boundedness. Here we want to remain as concrete as possible and look for practicable assumptions.

\mb We define operators $B_Q$ and Hardy spaces adapted to the operator $T$. Then using interpolation, we prove $L^p$ boundedness of this one. \\
 In particular case we will see that the $H^1_{F,\epsilon,mol}-L^1$ continuity of the operator $T$ below depends only on $L^2$ assumptions. It is only when we want to deduce $L^p$ estimates that we need stronger assumptions which imply $R$-boundedness used in \cite{weis}.

\mb Now we describe the choice of the collection ${\mathbb B}$, adapted to this operator. Then we will
check that the assumption (\ref{operh}) and the one about $M_{q_0}$
are satisfied. To finish the proof, we will show the $H^1_{F,\epsilon,
mol}-L^1$ boundedness of $T$ in Theorem \ref{theo2ap}.

\gb Equip $X=J \times Y$ with the parabolic quasi-distance $d$ and the measure $\mu$ defined by~:
$$d \Big( (t_1,y_1),(t_2,y_2) \Big) = \max\left\{ d_Y(y_1,y_2) , \sqrt{|t_1-t_2|} \right\} \qquad \textrm{and} \qquad  d\mu = dt \otimes d\nu.$$
If we write $\delta$ for the homogeneous dimension of the space $(Y,d_Y,\nu)$, then the space $X$ is of homogeneous type with homogeneous dimension $\delta+2$. We choose $\varphi \in \s(\R^+)$ such that $ \int_{ \R^+} \varphi(t)dt=1$ and $\varphi(t):=0$ for all $t<0$ ($\varphi$ does not need to be continuous at $0$). In fact we shall use only the fast decay of $\varphi$ and we will never consider regularity about it. In addition, we have added a condition for the support. This is a ``physical'' heuristics : this condition permits to define $A_Q(f)(t,x)$ by (\ref{opA}) with only $(f(\sigma,y))_{\sigma\leq t}$, which corresponds to the ``past informations'' about $f$. However we do not really need this assumption in the sequel. \\
For each ball $Q$ of $X$, we
write $r_Q$ its radius and we define the $B_Q$ operator as~:
\be{bq} B_Q=B_{r_Q^2} \qquad \textrm{with} \qquad B_r(f):= f-A_r(f), \ee
where the operator $A_r$ is defined by~: \be{opA} A_{r} (f)(t,x) :=\int _{\sigma =0}^{+\infty} \varphi_{r}(t-\sigma)
e^{rL}(f(\sigma,.))(x)  d\sigma. \ee Here we write $\varphi_{r}$ as the $L^1(\R)$ normalized function $\varphi_r(t):=r^{-1} \varphi(t/r)$. In fact, the integral for $\sigma\in[0,\infty)$ is reduced to $[0,t]$, due to the fact that $\varphi$ is supported in $\R^+$.

\mb Now to check the abstract assumption on the Hardy space, to be able to interpolate our operator $T$, we will use some conditions on our semigroup $e^{tL}$. We refer the reader to the work of P. Auscher and J.M. Martell (\cite{AM2}) to a precise study of off-diagonal estimates. Here we exactly define the decays, which will later be required. 

\begin{defn} \label{def:off} Let ${\mathcal T}:=(T_t)_{t\in J}$ be a collection of $L^2(Y)$-bounded operators and $p$ a positive integer. We will say that ${\mathcal T}$ satisfies  off-diagonal $L^2-L^2$ estimates at order $p$ if there exists a bounded function $\gamma$ satisfying
\be{gammaint} \forall\, 0\leq k\leq p, \qquad \sup_{u\geq 0} \gamma(u) (1+u)^{k} <\infty , \ee
such that for all balls $B\subset Y$ of radius $r$, for all functions $f$ supported in $B$ then
\be{offdiagonallow}
 \left( \frac{1}{\nu(2^{j+1}B)} \int_{S_j(B)} \left|T_{u^2} (f) \right|^{2} d\nu \right)^{1/2} \leq \frac{\nu(B)}{\nu(2^{j+1}B)} \gamma\left(\frac{2^{j+1}r}{u} \right) \left( \frac{1}{\nu(B)}\int_{B} |f|^2 d\nu \right)^{1/2}.\ee
We also write ${\mathcal T} \in {\mathcal O}_{p}(L^2-L^2)$. 
\end{defn}

\begin{rem}  This condition is satisfied for $p=\infty$ if the kernel $K_t$ of the operator $T_{t}$ admits some gaussian estimates like
$$ \left|K_t(x,y)\right| \lesssim \frac{1}{\nu(B(x,t^{1/2}))} e^{-\rho d(x,y)^2/t},$$
with $\rho>0$.
\end{rem}

\mb We will prove the following result in the next section~:

\begin{thm} \label{thresume} Let $L$ be a generator of a bounded analytic semigroup ${\mathcal T}:=(e^{tL})_{t>0}$ on $L^2(Y)$ such that
${\mathcal T},(tLe^{tL})_{t>0}$ and $(t^2L^2e^{tL})_{t>0}$ belong to ${\mathcal O}_{4}(L^2-L^2)$. Then for all $\epsilon>0$ the operator $T$ is continuous from $H^1_{F,\epsilon, mol}(X)$ to $L^1(X)$.
\end{thm}

\begin{rem} We recall that the semigroup $(e^{tL})_{t>0}$ is supposed to be analytic on $L^2$. Using Cauchy formula, if $(e^{zL})_{z}$ satisfies the $L^2-L^2$ off-diagonal estimates of ${\mathcal O}_{4}(L^2-L^2)$ for the complex variable $z$ belonging to a complex cone, then $(tLe^{tL})_{t>0}$ and $(t^2L^2e^{tL})_{t>0}$ belong to ${\mathcal O}_{4}(L^2-L^2)$.
\end{rem}

\mb We finish this Section by explaining how we can use this result to obtain positive answer for the maximal regularity problem. We want to apply the abstract results, recalled in the previous section. First we have to check the assumption (\ref{operh})~:

\begin{prop} There is a constant $0<A'<\infty$ so that for all $r>0$ the operator $A_r$ is $L^2(X)$ bounded and we have~:
$$ \left\| A_r \right\|_{L^2 \rightarrow L^2} \leq A'.$$
\end{prop}

\deme By definition the semigroup $e^{rL}$ is $L^2(Y)$-bounded so we have the following estimates~:
\begin{align*}
\left\| A_{r} (f) \right\|_{2} & \leq \left\| \int _{\sigma =0}^{+\infty} \int _{y \in Y} \left|\varphi_{r}(t-\sigma)\right|
\left\|e^{rL}(f(\sigma,.)) \right\|_{2,d\nu}  d\sigma \right\|_{2,dt} \\
 & \lesssim \left\| \int _{\sigma =0}^{+\infty}  \left|\varphi_{r}(t-\sigma)\right|
\left\|f(\sigma,.) \right\|_{2,d\nu}  d\sigma \right\|_{2,dt} \\
 & \lesssim \left\| \int _{\sigma =-\infty}^{t}  \left|\varphi_{r}(\sigma)\right|
\left\|f(t-\sigma,.) \right\|_{2,d\nu}  d\sigma \right\|_{2,dt} \\
 & \lesssim \left\|\varphi_{r}\right\|_1 \left\|f \right\|_{2,d\mu} \lesssim \left\|f \right\|_{2,d\mu}.
\end{align*}
So we have proved that $A_r$ is $L^2(X)$-bounded and its boundedness is uniform for $r>0$.
\findem

\begin{thm} The operator $T$ is $L^2(X)$-bounded.
\end{thm}

\mb This fact was proved in \cite{DS} because it is equivalent to the maximal $L^2$ regularity on $L^2(Y)$.

\mb Applying Theorem \ref{theogeneh}, we have~:

\begin{thm} Let $L$ be a generator of a bounded analytic semigroup ${\mathcal T}:=(e^{tL})_{t>0}$ on $L^2(Y)$ such that
${\mathcal T},(tLe^{tL})_{t>0}$ and $(t^2L^2e^{tL})_{t>0}$ belong to ${\mathcal O}_{4}(L^2-L^2)$. Let us assume that for $q_0\in(2,\infty]$
\be{hypi} M_{q_0} \lesssim M_{HL,2}. \ee
Then for all exponent $p\in(q_0',2]$, the operator $T$ admits a continuous extension on $L^p(Y)$ and so $L$ has maximal $L^p$-regularity.
In addition we can have weighted results~: let $\omega\in {\mathbb A}_\infty(Y)$ be a weight on $Y$. Then for all exponents $p\in {\mathcal W}_{\omega}(q_0',2)$
$L$ has maximal $L^p(\omega)$-regularity.
\end{thm}

\deme The first part of the theorem is a direct consequence of Theorem \ref{theogeneh} as the above assumptions were checked before. The second part about weighted results is an application of Theorem \ref{theopoids} with the following property. For $\omega\in {\mathbb A}_\infty(Y)$ a weight on the space $Y$, we set $\tilde{\omega}$ for the associated weight on $X=J \times Y$ defined by the tensor product $\tilde{\omega}:={\bf 1_R} \otimes \omega$ ~: for all ball $Q\subset X$ of radius $r_Q$, we can write $Q=I\times Q_Y$ with an interval $I$ of length $r_Q^2$ and $Q_Y\subset Y$ a ball of radius $r_Q$ 
$$ \tilde{\omega}(Q):=r_Q^2 \omega(Q_Y).$$
Then with this definition, it is obvious to check that for exponents $p,q\in[1,\infty]$~:
$$ \tilde{\omega}\in {\mathbb A}_p(X)  \Longleftrightarrow \omega \in {\mathbb A}_p(Y) $$
and
$$\tilde{\omega}\in RH_q(X)  \Longleftrightarrow \omega \in RH_q(Y). $$
So we have 
$$ \mathcal{W}_{\omega}(\sigma',2)= \mathcal{W}_{\tilde{\omega}}(\sigma',2).$$
\findem

\mb We want now to study the main assumption (\ref{hypi}). For example, we give  other stronger assumption, describing ``off-diagonal'' estimates.

\begin{prop} \label{assumcas} We recall the maximal operator
$$ M_{q_0} (f)(\sigma,x):= \sup_{\genfrac{}{}{0pt}{}{Q \textrm{ ball}}{(\sigma,x)\in Q}} \left( \frac{1}{\mu(Q)} \int_{Q} \left| A_Q^*(f)\right|^{q_0} d\mu \right)^{1/q_0}.$$
If the semigroup $(e^{tL^*})_{t>0}$ satisfy these ``$L^{q_0'}-L^2$ off-diagonal'' estimates~:
there exist coefficients $(\beta_j)_{j\geq 0}$ satisfying
\be{beta} \sum_{j\geq 0} 2^{j} \beta_j <\infty \ee
such that for all balls $B$ and for all functions $f\in L^2(Y)$ we have
 \be{offdiag} \left( \frac{1}{\nu(B)} \int_B \left|e^{r_B^2L^*} (f) \right|^{q_0'} d\nu \right)^{1/q_0'} \leq \sum_{j\geq 0} \beta_{j} \left( \frac{1}{\nu(2^{j} B)}\int_{2^{j}B} |f|^2 d\nu \right)^{1/2}.\ee
Then $M_{q_0}$ is bounded by the Hardy-Littlewood maximal operator $M_{HL,2}$ on $X$, so (\ref{hypi}) is satisfied.
\end{prop}

\deme  Let $Q$ be a ball containing the point $(\sigma,x)\in X$ and $r_Q$ be its radius.
For $f,g\in L^2(X)$ we have~:
\begin{align*}
\langle A_{Q} (f) ,g \rangle & := \int_{(t,x)\in X} \int _{\sigma =0}^{+\infty} \varphi_{r_Q^2}(t-\sigma)
e^{r_Q^2 L}(f(\sigma,.))(x) g(t,x)  d\sigma dt d\nu(x) \\
 & = \int_{(t,x)\in X} \int _{\sigma =0}^{+\infty} \varphi_{r_Q^2}(t-\sigma)
f(\sigma,x) \left[\left(e^{r_Q^2 L}\right)^*g(t,.)\right](x)  d\sigma dt d\nu(x).
\end{align*}
So we conclude that~: \be{Aetoile} A_Q^*(g)(\sigma,x):= \int_{t\in
\R^+}  \varphi_{r_Q^2}(t-\sigma) \left[\left(e^{r_Q^2
L}\right)^*g(t,.)\right](x) dt.\ee By using the Minkowski
inequality, we also have  that
$$ \left ( \int_{Q} \left| A_Q^*(g)\right|^{q_0} d\mu \right)^{\frac{1}{q_0}}  \leq \int_{t\in \R^+}  \left\| \varphi_{r_Q^2}(t-\sigma) \left[\left(e^{r_Q^2 L}\right)^*g(t,.)\right](x) {\bf 1}_{Q}(\sigma,x) \right\|_{q_0,d\nu(x)d\sigma} dt. $$
By definition of the parabolic quasi-distance, we can write
$$ Q = I \times B$$
with $I$ an interval of lenght $r_Q^2$ and $B$ a ball of $Y$ of radius $r_Q$. Then we have~:
\begin{align*}
\lefteqn{ \left( \int_{Q} \left| A_Q^*(f)\right|^{q_0} d\mu \right)^{1/q_0}  \leq} & & \\
& &  \int_{t\in \R^+}  \left\| \varphi_{r_Q^2}(t-\sigma){\bf 1}_{I}(\sigma) \right\|_{q_0,d\sigma} \left\|{\bf 1}_{B}(x) \left(e^{r_Q^2 L}\right)^*g(t,.)(x) \right\|_{q_0,d\nu(x)} dt.
\end{align*}
With the assumption (\ref{offdiag}), we obtain
\begin{align*}
\lefteqn{ \left( \int_{Q} \left| A_Q^*(f)\right|^{q_0} d\mu \right)^{1/q_0} \leq } & & \\
 & &   \sum_{j\geq 0} \int_{t\in \R^+}  \left\| \varphi_{r_Q^2}(t-\sigma) {\bf 1}_{I}(\sigma)\right\|_{q_0,d\sigma}  \beta_j \frac{\nu(B)^{1/q_0}}{\nu(2^{j}B)^{1/2}} \left\|g(t,x) {\bf 1}_{2^{j}B}(x)\right\|_{2,d\nu(x)} dt.
 \end{align*}
Now we decompose the integration over $t$ by~:
\begin{align*}
\lefteqn{ \left( \int_{Q} \left| A_Q^*(f)\right|^{q_0} d\mu \right)^{1/q_0}  \leq} & & \\
 & & \sum_{j\geq 0} \sum_{k\geq 0} \int_{t \in S_k(I)}  \left\| \varphi_{r_Q^2}(t-\sigma) {\bf 1}_{I}(\sigma) \right\|_{q_0,d\sigma}  \beta_j \frac{\nu(B)^{1/q_0}}{\nu(2^{j}B)^{1/2}} \left\|g(t,x) {\bf 1}_{2^jB}(x) \right\|_{2,d\nu(x)} dt.
 \end{align*}
With the Cauchy-Schwarz inequality, we have
\begin{align*}
\lefteqn{\left( \int_{Q} \left| A_Q^*(f)\right|^{q_0} d\mu \right)^{1/q_0} } & & \\
& & \lesssim \sum_{j\geq 0} \sum_{k\geq 0} r_Q^{-2}
\left(1+ 2^{k}\right)^{-l} r_Q^{2/q_0} \beta_j \frac{\nu(B)^{1/q_0}}{\nu(2^{j}B)^{1/2}} (2^{k} r_Q^2)^{1/2}\left\|g(t,x){\bf 1}_{2^kI \times 2^jB}(t,x) \right\|_{2,dt d\nu(x)} \\
& & \lesssim \sum_{j\geq 0} \sum_{k\geq 0} r_Q^{-1+2/q_0}
\left(1+ 2^{k}\right)^{-l+1/2} \beta_j \frac{\nu(B)^{1/q_0}}{\nu(2^{j}B)^{1/2}} \left\|g(t,x) {\bf 1}_{2^kI \times 2^jB}(t,x)\right\|_{2,dtd\nu(x)}.
\end{align*}
Here $l$ is an integer as large as we want, due to the fast decay of $\varphi$. Using the  Hardy-Littlewood maximal operator, we have
$$ \left\|g(t,x) {\bf 1}_{2^kI \times 2^jB}(t,x) \right\|_{2,dtd\nu(x)} \leq \mu\left(\max\{2^j,2^{k/2}\} Q\right)^{1/2}  \inf_{Q} M_{HL,2}(g).$$
So we obtain
\begin{align*}
\lefteqn{\left( \int_{Q} \left| A_Q^*(g)\right|^{q_0} d\mu \right)^{1/q_0} \leq} & & \\
& &  \left[ \sum_{j\geq 0} \sum_{k\geq 0} r_Q^{-1+2/q_0} \left(1+ 2^{k}\right)^{-l+1/2} \beta_j \frac{\nu(B)^{1/q_0}}{\nu(2^{j}B)^{1/2}}  \mu\left(\max\{2^j,2^{k/2}\} Q\right)^{1/2} \right]  \inf_{Q} M_{HL,2}(g). \label{equ1}
\end{align*}
We now estimate the sum over the parameters $j$ and $k$. We have the two following cases.
Write
$$  S_1 := \sum_{j\geq k/2 \geq 0} r_Q^{-1+2/q_0} \left(1+ 2^{k}\right)^{-l+1/2} \beta_j \frac{\nu(B)^{1/q_0}}{\nu(2^{j}B)^{1/2}}  \mu\left( 2^{j} Q\right)^{1/2} $$
and
$$ S_2 := \sum_{k/2 \geq j \geq 0} r_Q^{-1+2/q_0} \left(1+ 2^{k}\right)^{-l+1/2} \beta_j \frac{\nu(B)^{1/q_0}}{\nu(2^{j}B)^{1/2}}  \mu\left(2^{k/2} Q\right)^{1/2} .$$
We must estimate these two sums. For the first, we use that $\mu(Q)= |I| \nu(B)=r_Q^2 \nu(B)$ to have
\begin{align*}
  S_1 & \leq  \sum_{j\geq k/2 \geq 0} 2^{j} \left(1+ 2^{k}\right)^{-l+1/2} \beta_j \frac{\mu(Q)^{1/q_0}}{\mu(2^{j}Q)^{1/2}}  \mu\left( 2^{j} Q\right)^{1/2} \\
   & \leq  \mu(Q)^{1/q_0} \sum_{j\geq k/2 \geq 0} 2^{j} \left(1+ 2^{k}\right)^{-l+1/2} \beta_j \\
   & \leq  \mu(Q)^{1/q_0} \sum_{j\geq 0} 2^{j} \beta_j \lesssim \mu(Q)^{1/q_0}. \\
\end{align*}
In the last inequality, we have used the assumption (\ref{beta}) about the coefficients $(\beta_j)_j$. \\
For the second sum, we have (with the doubling property of $\mu$ and $l$ large enough)
 \begin{align*}
  S_2 & \leq r_Q^{2/q_0} \nu(B)^{1/q_0} \sum_{k/2 \geq j \geq 0} r_Q^{-1} \left(1+ 2^{k}\right)^{-l+1/2} \beta_j \left(\frac{ \mu(2^{k/2} Q) }{\nu(2^{j}B)}\right)^{1/2}  \\
   & \lesssim \mu(Q)^{1/q_0} \sum_{k/2 \geq j \geq 0}  r_Q^{-1} \left(1+ 2^{k}\right)^{-l+1/2} \beta_j \left(\frac{ \mu(2^{j} Q) }{\nu(2^{j}B)}\right)^{1/2} 2^{(k/2-j)(\delta+2)/2} \\
   & \lesssim \mu(Q)^{1/q_0} \sum_{k/2 \geq j \geq 0}  \left(1+ 2^{k}\right)^{-l+1/2} \beta_j 2^{j} 2^{(k/2-j)(\delta+2)/2}  \\
   & \lesssim \mu(Q)^{1/q_0} \sum_{j \geq 0}  \left(1+ 2^{j}\right)^{-l+4+\delta/2} \beta_j 2^{-j(\delta/2+1)} \lesssim \mu(Q)^{1/q_0}.
\end{align*}
So we have proved that there exists a constant $C$ (independant on $g$ and $Q$) such that~:
$$ \left( \int_{Q} \left| A_Q^*(g)\right|^{q_0} d\mu \right)^{1/q_0} \leq C \mu(Q)^{1/q_0}  \inf_{Q} M_{HL,2}(g). $$
 We can also conclude that
 $$M_{q_0}(f) \lesssim M_{HL,2}(g).$$
  \findem

\mb We have described an ``off-diagonal'' estimates implying (\ref{hypi}) with dyadic scale. Obviously we could describe other ``off-diagonal'' estimates ... 

\mb We would like to finish this section by comparing our result with the one of S. Blunck and P.C. Kunstmann \cite{BK2}. In their paper, the authors have used their assumptions (\ref{assum1}) and (\ref{assum2}) to use inside their proof the following inequality~: for all ball $Q_Y$ of $Y$ and all functions $f\in L^2(Y)$
\be{BK} \left(\frac{1}{\nu(Q_Y)}\int_{Q_Y} \left|e^{-r_{Q_Y}^2L}f\right|^{2} d\nu \right)^{1/2} \lesssim \left(\sum_{k=0}^\infty k^{-1-\epsilon} \left[\frac{1}{\nu(kQ_Y)}\int_{kQ_Y} |f|^{p_0} d\nu\right] \right)^{1/p_0}.\ee
With this inequality, a simple computation gives us that for all function $f\in L^2(X)$ and all balls $Q$ of $X$
 \be{hypii} \left(\frac{1}{\nu(Q)}\int_{Q} \left|A_Q(f)\right|^{2} d\mu \right)^{1/2} \lesssim \inf_{Q} M_{HL,p_0}(f). \ee
It is surprizing to note that their assumption (\ref{hypii}) seems to be not comparable with ours (\ref{hypi}). These two assumptions are quite different in the sense that we require different kind of ``off-diagonal'' estimates, however they seem to be the dual of each other.

\section{Proof of Theorem \ref{thresume}.} \label{section4}

 This section is devoted to the proof of a technical result : Theorem \ref{thresume}. Let us first repeat it.

\begin{thm} \label{theo2ap}  Let $L$ be a generator of a bounded analytic semigroup on $L^2(Y)$. Assume that $(e^{tL})_{t>0},(tLe^{tL})_{t>0}$ and $(t^2L^2e^{tL})_{t>0}$ belong to ${\mathcal O}_{4}(L^2-L^2)$. Then there exist coefficients $\alpha_{j,k}$ such that for all balls $Q\subset X$, for all $k\geq 0,j\geq 2$ and for all functions $f$ supported in $S_k(Q)$
\be{theo2a} \left(\frac{1}{\mu(2^{j+k+1}Q)} \int_{S_j(2^kQ)} \left| T(B_Q(f))\right|^2 d\mu \right)^{1/2} \leq \alpha_{j,k} \left(\frac{1}{\mu(2^{k+1}Q)}\int_{S_k(Q)} |f|^2 d\mu \right)^{1/2}. \ee
In addition the coefficients $\alpha_{j,k}$ (independent in $Q$) satisfy \be{hypap} \Lambda:= \sup_{Q} \ \sup_{k\geq 0}
\left[\sum_{j\geq 2} \frac{\mu(2^{j+k+1}Q)}{\mu(2^{k+1}Q)} \alpha_{j,k} \right] <\infty. \ee
With Theorem \ref{theo2h}, these estimates imply the $H^1_{F,\epsilon, mol}(X)-L^1(X)$ boundedness of $T$ for all $\epsilon>0$.
\end{thm}

\deme We write $r=r_Q$ and $(t_0,x_0)$ the radius and the center of the ball $Q$ so we have defined $B_Q$ as $B_{r^2}$. The function $f$ is fixed. The parameter $j$ and $k$ are fixed too. We write $Q$ as the product $Q=I \times B$ with $I$ an interval of length $r_Q^2$ and $B$ a ball of $Y$ of radius $r_Q$.
We have
\begin{align*}
T B_{r^2}(f) (t,x) & = T(f)(t,x) - T A_{r^2}(f) (t,x ) \\
 & =\int _{0} ^{t} \left[Le^{(t-s)L}f(s,.)\right](x) ds-\int_{0}^{t} \left[L e^{(t-s)L} A_{r^2} f(s,.)\right](x) ds,
\end{align*}
where
\begin{align*}
 \left[L e^{(t-s)L} A_{r^2}f(s,.) \right](x) = L e^{(t-s)L} \left[\int_{\sigma=0}^{+\infty} \varphi
_{r^2} (s-\sigma)  e^{r^2 L} f(\sigma,.)  d\sigma \right](x).
\end{align*}
So we obtain
\begin{align}
\lefteqn{T(B_{r^2}f)(t,x)  =} & & \nonumber \\
 & & \int_\R \int_\R \varphi_{r^2}(s-\sigma)\left[{\bf 1}_{0<\sigma\leq t} Le^{(t-\sigma)L}f(\sigma,.)(x) -{\bf 1}_{0<s\leq t} L e^{(t-s+r^2)L}f(\sigma,.)(x) \right] d\sigma ds . \label{importante}
 \end{align}
We have three time-parameters $\sigma,t$ and $s$. As in the case of Calder\'on-Zygmund operators, the difference within the two brackets is very important. This will allow us to obtain the necessary decay for the coefficients $\alpha_{j,k}$.
We decompose into two domains~:
$$ D_1:=\left\{ (\sigma,t,s),\ 0\leq \sigma \leq t\leq s \right\} \quad \textrm{and} \quad D_2:=\left\{ (\sigma,t,s),\ 0\leq s,\sigma \leq t \right\}.$$
For $i\in\{1,2\}$ we set $D_i(t):=\left\{(\sigma,s);\ (\sigma,t,s)\in D_i\right\}$ and
\begin{align*}
\lefteqn{ U_i(f)(t,x) :=} & & \\
 & &  \iint_{D_i(t)} \varphi_{r^2}(s-\sigma)\left[{\bf 1}_{0<\sigma\leq t} Le^{(t-\sigma)L}f(\sigma,.)(x) -{\bf 1}_{0<s\leq t} L e^{(t-s+r^2)L}f(\sigma,.)(x) \right] d\sigma ds.
 \end{align*}
As $\varphi$ is supported in $\R^+$, we have decomposed
\be{decompositionTB} T(B_{r^2}f)(t,x) = \sum_{i=1}^{2} U_i(f)(t,x). \ee
If we do not want to use the condition of the support of $\varphi$, there is a third term which is estimated as the first one. \\
We begin the study when one of the two terms, in the square brackets, vanishes. The radius $r$ is fixed for all the proof and we set
$$ \chi_N(y):=\frac{1}{r^2}\left(1+ \frac{|y|}{r^2} \right)^{-N}.$$
 \\
$1-)$ First case : $(\sigma,s)\in D_1(t)$. \\
Here we have the following expression~:
\begin{align*}
U_1(f)(t,x) = \int_{t}^\infty \int_{0}^{t} \varphi_{r^2}(s-\sigma) L e^{(t-\sigma)L} f(\sigma,.)(x) d\sigma ds.
\end{align*}
There is no ``cancellation'' so we can directly estimate it by using the fast decay of $\varphi$.
For $N$ a large enough integer
\begin{align*}
\left| U_1(f)(t,x) \right| & \lesssim  r^2\int_{0}^{t} \chi_N(t-\sigma) \left| L e^{(t-\sigma)L} f(\sigma,.)(x)\right| d\sigma.
\end{align*}
By definition of the parabolic quasi-distance,
$$(t,x) \in S_j (2^k Q) \Longleftrightarrow  \left\{ \begin{array}{l}
d_Y(x,x_0) \simeq 2^{k+j}r \\
|t-t_0| \leq ( 2^{k+j+1}  r)^2 \end{array} \right.
\textrm{ or }
\left\{\begin{array}{l}
d_Y(x,x_0) \lesssim 2^{k+j}r \\
|t-t_0| \simeq (2^{k+j+1} r)^2 \end{array} \right. .
$$
So, as $f$ is supported in $2^kQ$, we have
\be{III} \left(\frac{1}{\mu(2^{j+k+1}Q) }\int_{S_j(2^kQ)} \left| U_1(f)(t,x) \right|^2 d\nu(x)dt \right)^{1/2} \lesssim I+ II \ee
with
$$ I := \left( \frac{r^4}{\mu(2^{j+k+1}Q)} \int_{2^{2(j+k)}I} \left( \int_{2^{2k}I} \chi_N(t-\sigma) \left\| L e^{(t-\sigma)L} f(\sigma,.)\right\|_{2,S_j(2^kB)} d\sigma\right)^2 dt \right)^{1/2}$$
and
$$ II :=  \left(\frac{1}{\mu(2^{j+k+1}Q)} \int_{S_{2j}(2^{2k}I)} \left(\int_{2^{2k}I} \frac{ \left\| L e^{(t-\sigma)L} f(\sigma,.)\right\|_{2,2^{k+j+1}B}}{2^{2N(k+j)}} d\sigma\right)^2 dt \right)^{1/2}.$$

\noindent $*$ Study of $I$. \\
By using off-diagonal estimates $L^2-L^2$ (\ref{offdiagonallow}), we know that
\begin{align*}
 \lefteqn{\frac{1}{\nu(2^{j+k+1}B)^{1/2}} \left\| L e^{(t-\sigma)L} f(\sigma,.)\right\|_{2,S_j(2^kB)} \leq} & & \\
 & &  \frac{\nu(2^{k+1}B)}{\nu(2^{k+j+1}B)|t-\sigma|} \gamma\left(\frac{2^{j+k}r}{\sqrt{|t-\sigma|}}\right)  \left( \frac{1}{\nu(2^{k+1}B)} \int_{2^{k+1} B} |f(\sigma,.)|^2 d\nu\right)^{1/2}.
 \end{align*}
That is why, by using Cauchy-Schwarz inequality and the equality
$$ \mu(2^{j+k+1}Q) = \nu(2^{j+k+1}B) 2^{2(j+k)}r^2,$$
we estimate $I$ by the product
\begin{align*}
  \left( \frac{1}{2^{2(k+j)}} \int_{2^{2k+2j}I} \int_{2^{2k}I} \chi_{2N}(t-\sigma) \left(\frac{\nu(2^{k+1}B)}{\nu(2^{k+j+1}B)|t-\sigma|}\right)^2 \gamma\left(\frac{2^{j+k}r}{\sqrt{|t-\sigma|}}\right)^2  d\sigma dt \right)^{1/2}      \\
    \qquad  2^{k}r\left( \frac{1}{\mu(2^{k+1}Q)} \int_{2^{k+1} Q} |f|^2 d\mu\right)^{1/2}.
\end{align*}
Then we get
\begin{align*}
  I & \lesssim \frac{\nu(2^{k+1}B)}{\nu(2^{k+j+1}B)} \Bigg[ \frac{1}{|2^{2k+2j} I|} \int_{2^{2k+2j}I} \int_{2^{2k}I} \chi_{2N}(t-\sigma)  \frac{2^{2k}r^4}{|t-\sigma|^2}  \\
    & \hspace{2cm}  \gamma\left(\frac{2^{j+k}r}{\sqrt{|t-\sigma|}}\right)^2 d\sigma dt \Bigg]^{1/2} \left( \frac{1}{\mu(2^{k+1}Q)} \int_{2^{k+1} Q} |f|^2 d\mu\right)^{1/2} \\
  & \lesssim \frac{2^k \nu(2^{k+1}B)}{\nu(2^{k+j+1}B)} \left(\int_{0}^{2^{2(j+k)}} \left(1+v \right)^{-2N}  \frac{1}{v^2} \gamma\left(\frac{2^{j+k}}{\sqrt{v}}\right)^2 dv \right)^{1/2} \\
  &  \hspace{5cm} \left( \frac{1}{\mu(2^{k+1}Q)} \int_{2^{k+1} Q} |f|^2 d\mu\right)^{1/2} \\
  & \lesssim \frac{2^{k}\nu(2^{k+1}B)}{\nu(2^{k+j+1}B)} 2^{-j-k}\left(\int_{1}^{\infty} \left(1+2^{k+j}v^{-2} \right)^{-2N} \gamma(v)^2 v dv\right)^{1/2} \\
  &  \hspace{5cm} \left( \frac{1}{\mu(2^{k+1}Q)} \int_{2^{k+1} Q} |f|^2 d\mu\right)^{1/2}.
\end{align*}
$*$ Study of $II$. \\
In this case, we have $t\in S_{2j}(2^{2k}I)$ and $\sigma \in 2^{2k}I$, so
 $$ |t-\sigma|\simeq 2^{2(j+k)}r^2.$$
By using off-diagonal estimates (\ref{offdiagonallow}), we know that
\begin{align*}
 \lefteqn{\frac{1}{\nu(2^{j+k+1}B)^{1/2}} \left\| L e^{(t-\sigma)L} f(\sigma,.)\right\|_{2,2^{k+j+1}B} \lesssim \hspace{2cm} } & & \\
 & &  \hspace{2cm} \frac{1}{2^{2(j+k)} r^2}  \left( \frac{1}{\nu(2^{k+1}B)} \int_{2^{k+1} B} |f(\sigma,.)|^2 d\nu\right)^{1/2}.
 \end{align*}
So we obtain that
\begin{align*}
  II & \lesssim   \left(\frac{1}{2^{2(j+k)}r^2} \int_{t\in S_{2j}(2^{2k}I)} \int_{2^{2k}I} 2^{-4N(k+j)} \frac{2^{2k}r^2 }{2^{4(j+k)} r^4}  d\sigma dt \right)^{1/2} & \\
  & & \hspace{-3cm} \left( \frac{1}{\mu(2^{k+1}Q)} \int_{2^{k+1} Q} |f|^2 d\mu\right)^{1/2} \\
  & \lesssim  2^{-2j}2^{-2N(k+j)} \left( \frac{1}{\mu(2^{k+1}Q)} \int_{2^{k+1} Q} |f|^2 d\mu\right)^{1/2}.
\end{align*}

\gb We have also the following estimate
\begin{align*} \lefteqn{I+ II \lesssim \left( \frac{1}{\mu(2^{k+1}Q)} \int_{2^{k+1} Q} |f|^2 d\mu\right)^{1/2} } & & \\
 & &    \left[2^{-2j} 2^{-2N(k+j)}  + \frac{\nu(2^{k+1}B)}{\nu(2^{k+j+1}B)} 2^{-j-k} \left(\int_{1}^\infty \left(1+2^{k+j}v^{-2} \right)^{-2N} \gamma(v)^2 v dv \right)^{1/2}\right].
 \end{align*}
With (\ref{III}), here we can choose
$$ \alpha_{j,k}= \left[2^{-N(k+j)}  +  \frac{\nu(2^{k+1}B)}{\nu(2^{k+j+1}B)}2^{-j-k}\left(\int_{1}^\infty \left(1+2^{j}v^{-2} \right)^{-N} \gamma(v)^2 vdv \right)^{1/2} \right]$$
for $N$ a large enough integer. \\
$2-)$ Last case for $(\sigma,s)\in D_2(t)$ : $0\leq s,\sigma \leq t$.\\
The relation (\ref{importante}) gives us that~:
$$ U_2(f)(t,x)  =
\int_0^t \int_0^t \varphi_{r^2}(s-\sigma)\left[Le^{(t-\sigma)L}f(\sigma,.)(x) - L e^{(t-s+r^2)L}f(\sigma,.)(x) \right] d\sigma ds. $$
Here we use the time regularity. We have~:
\begin{align*}
\left|Le^{(t-\sigma)L}f(\sigma,.)(x) - L e^{(t-s+r^2)L}f(\sigma,.)(x) \right| & = \left| \int_{t-s+r^2}^{t-\sigma} \frac{\partial L e^{zL}f(\sigma,.)(x)}{\partial z} dz \right| \\
 & = \left| \int_{t-s+r^2}^{t-\sigma} L^2 e^{zL}f(\sigma,.)(x) dz \right|.
\end{align*}
We do not know whether $t-s+r^2 \leq t-\sigma$ or $t-s+r^2 \geq t-\sigma$, fortunately this order is not important.
Then we repeat the same arguments as before~: \\
\be{IIII} \left(\frac{1}{\mu(2^{j+k+1}Q)}\int_{S_j(2^kQ)} \left| U_2(f)(t,x) \right|^2 d\nu(x)dt \right)^{1/2} \lesssim I+ II \ee
with
\begin{align*}
 \lefteqn{I := \Bigg( \frac{1}{\mu(2^{j+k+1}Q)} \int_{2^{2k+2j}I} \left( \int_{2^{2k}I} \int_0^t \chi_N(s-\sigma) \right.  \hspace{3cm} } & & \\
 & & \left. \left. \hspace{3cm} \int_{t-s+r^2}^{t-\sigma} \left\| L^2 e^{zL} f(\sigma,.)\right\|_{2,S_j(2^kB)} dz ds d\sigma\right)^2 dt \right)^{1/2}
 \end{align*}
and
\begin{align*}
\lefteqn{ II :=  \left( \frac{1}{\mu(2^{j+k+1}Q)} \int_{t\in S_{2j}(2^{2k}I)} \left(\int_{2^{2k}I} \int_0^t \chi_N(s-\sigma)  \right. \right. \hspace{3cm} } & & \\
 & & \left. \left. \hspace{3cm} \int_{t-s+r^2}^{t-\sigma} \left\| L^2 e^{zL} f(\sigma,.)\right\|_{2,2^{k+j}B} dz ds d\sigma\right)^2 dt \right)^{1/2}.
 \end{align*}
$*$ Study of $I$. \\
By using off-diagonal estimates (\ref{offdiagonallow}), we know that
\begin{align*}
 \lefteqn{\frac{1}{\nu(2^{j+k+1}B)^{1/2}}\left\| L^2 e^{zL} f(\sigma,.)\right\|_{2,S_j(2^kB)} \leq \hspace{2cm} } & & \\
 & &  \hspace{2cm} \frac{\nu(2^kB)}{\nu(2^{k+j+1}B)z^2} \gamma\left(\frac{2^{j+k}r}{\sqrt{z}} \right)  \left( \frac{1}{\nu(2^{k}B)} \int_{2^k B} |f(\sigma,.)|^2 d\nu\right)^{1/2}.
\end{align*}
So we obtain
\begin{align*}
 \lefteqn{ I \lesssim  \frac{\nu(2^kB)}{\nu(2^{k+j+1}B)} \left( \frac{1}{\mu(2^{k}Q)} \int_{2^k Q} |f|^2 d\mu\right)^{1/2} {\Bigg [} \frac{1}{2^{2(k+j)} r^2} \int_{2^{2k+2j}I}\int_{2^{2k}I}  } & & \\
  & &  \left( \int_{0}^{t}  \frac{1}{r^2} \left(1+\frac{|s-\sigma|}{r^2}\right)^{-N} \left|\int_{t-s+r^2}^{t-\sigma}  \frac{2^{k}r}{z^2} \gamma\left(\frac{2^{j+k}r}{\sqrt{z}}\right) dz\right| ds\right)^2 d\sigma dt {\Bigg ]} ^{1/2} .
 \end{align*}
We use the fast decays of $\gamma$ and a large enough exponent $p\geq 4$ to obtain the following inequality
\begin{align*}
 \left| \int_{t-s+r^2}^{t-\sigma}  \frac{1}{z^2} \gamma\left(\frac{2^{j+k}r}{\sqrt{z}}\right) dz \right| & \lesssim  \left| \int_{t-s+r^2}^{t-\sigma}  \frac{1}{z^2} \left(\frac{2^{j+k}r}{\sqrt{z}}\right)^{-p} dz \right| \lesssim \left| \int_{t-s+r^2}^{t-\sigma}  \frac{z^{p/2-2}}{2^{p(j+k)}r^p} dz \right| \\
 & \lesssim 2^{-p(j+k)}r^{-p} \left|\sigma-s+r^2\right| \left(2^{j+k}r\right)^{2(p/2-2)}.
\end{align*}
At the last inequality, we have used that $|t-\sigma|\lesssim \left(2^{j+k}r\right)^{2}$ and similarly for $|\sigma-s|$.
So for example taking $p=4$, we obtain
\begin{align}
 \left| \int_{t-s+r^2}^{t-\sigma}  \frac{1}{z^2} \gamma\left(\frac{2^{j+k}r}{\sqrt{z}}\right) dz \right| & \lesssim 2^{-4(j+k)} r^{-2}\left(1+ \frac{|s-\sigma|}{r^2} \right). \label{difference}
 \end{align}
We also get
\begin{align*}
 I & \lesssim \frac{\nu(2^kB)2^{k}}{\nu(2^{k+j+1}B)2^{4(k+j)}r}  \left( \frac{1}{|2^{2k+2j} I|} \int_{2^{2k+2j}I}  2^{2k}r^2 dt \right)^{1/2} \left( \frac{1}{\mu(2^{k}Q)} \int_{2^k Q} |f|^2 d\mu\right)^{1/2} \\
  & \lesssim \frac{\nu(2^kB)}{\nu(2^{k+j+1}B)} 2^{-4j} \left( \frac{1}{\mu(2^{k}Q)} \int_{2^k Q} |f|^2 d\mu\right)^{1/2}.
\end{align*}
Here we can choose
$$ \alpha_{j,k}= \frac{\nu(2^kB)}{\nu(2^{k+j+1}B)} 2^{-4j}  \simeq \frac{\mu(2^kQ)}{\mu(2^{k+j+1}Q)} 2^{-2j} \|\gamma\|_\infty .$$
$*$ Study of $II$. \\
In this case, we have $|t-\sigma| \simeq 2^{2(k+j)}r^2$.
By using off-diagonal estimates (\ref{offdiagonallow}), we know that for $z\geq r^2$
$$ \frac{1}{\nu(2^{j+k+1}B)^{1/2}}\left\| L^2 e^{zL} f(\sigma,.)\right\|_{2,S_{j}(2^{k}B)} \leq \frac{1}{z^2}
  \frac{\nu(2^{k}B)}{\nu(2^{j+k+1}B)} \left( \frac{1}{\nu(2^{k}B)} \int_{2^k B} |f(\sigma,.)|^2 d\nu\right)^{1/2}. $$
So we obtain
\begin{align*} II & \lesssim  \frac{\nu(2^{k}B)}{2^{j} \nu(2^{j+k+1}B)} \left(\int_{2^{2k+2j}I} \int_{2^{2k}I} \left(\int_{0}^{t}  \chi_N(s-\sigma) \int_{t-s+r^2}^{t-\sigma}  \frac{1}{z^2} dz ds \right)^2 d\sigma dt \right)^{1/2}  \\
  & \hspace{7cm}  \left( \frac{1}{\mu(2^{k}Q)} \int_{2^k Q} |f|^2 d\mu\right)^{1/2}.
\end{align*}
Then we use
\begin{align*}
 \left| \int_{t-s+r^2}^{t-\sigma}  \frac{1}{z^2} dz \right| & \lesssim \frac{r^2 + |s-\sigma|}{(t-\sigma)(t-s+r^2)}
\lesssim \frac{1 + \frac{|s-\sigma|}{r^2}}{2^{2(k+j)}r^2(1+\frac{|t-s|}{r^2})} \\
 & \lesssim \frac{\left(1+\frac{|s-\sigma|}{r^2}\right)^2}{2^{2(k+j)}r^2(1+\frac{|t-\sigma|}{r^2})} \lesssim \frac{\left(1+\frac{|s-\sigma|}{r^2}\right)^2}{2^{4(k+j)}r^2}
\end{align*}
 to finally obtain (with an other exponent $N$)
\begin{align*} II & \lesssim  \frac{\nu(2^{k}B)}{\nu( 2^{j+k+1}B)2^{j} 2^{4(j+k)} r^2} \left( \int_{2^{2k+2j}I} \int_{2^{2k}I} \left( \int_{0}^{t} \chi_N(s-\sigma) ds \right)^2 d\sigma dt \right)^{1/2}  \\
  & \hspace{7cm}  \left( \frac{1}{\mu(2^{k}Q)} \int_{2^k Q} |f|^2 d\mu\right)^{1/2} \\
  & \lesssim \frac{\nu(2^{k}B)2^{2k} r^2}{\nu( 2^{j+k+1}B)2^{4(j+k)} r^2}\left( \frac{1}{\mu(2^{k}Q)} \int_{2^k Q} |f|^2 d\mu\right)^{1/2} \\
  & \lesssim \frac{\mu(2^{k}Q)}{\mu( 2^{j+k+1}Q)2^{2j+2k}}\left( \frac{1}{\mu(2^{k}Q)} \int_{2^k Q} |f|^2 d\mu\right)^{1/2}.
\end{align*}
So here we can choose
$$ \alpha_{j,k}= \frac{\mu(2^{k}Q)}{\mu( 2^{j+k+1}Q)2^{2j+2k}}. $$

\noindent $3-)$ End of the proof.\\
With the decomposition (\ref{decompositionTB}), we have proved in the two previous points that we have the estimate
(\ref{theo2a}) with the coefficients $\alpha_{j,k}$ satisfying
\begin{align*}
 \alpha_{j,k} \lesssim & \  2^{-N(k+j)}  + \frac{\nu(2^{k}B)}{\nu(2^{k+j+1}B)}2^{-j} \left(\int_{1}^\infty \left(1+2^{j}v^{-2} \right)^{-2N}\gamma(v)^2 vdv \right)^{1/2}  \\
 &  \ + \frac{\mu(2^{k}Q)}{\mu( 2^{j+k+1}Q)2^{2j}}.
\end{align*}
We are going to check that (\ref{hypap}) is satisfied. So we must bound the quantity
$$ \lambda_{k,Q}:= \sum_{j\geq 2} \frac{\mu(2^{j+k+1}Q)}{\mu(2^{k+1}Q)} \alpha_{j,k}  $$
by a constant (independent on $k$ and $Q$).
The coefficient $\alpha_{j,k}$ is estimated by three terms. By using the doubling property for $\mu$, with $N$ large enough we can sum the first term $2^{-N(k+j)}$. For the second term with $N\geq 2$, we use (\ref{gammaint}) to have
\begin{align*}
\sum_{j\geq 2} \frac{\mu(2^{j+k+1}Q)}{\mu(2^{k+1}Q)} \frac{\nu(2^{k}B)}{\nu(2^{k+j+1}B)}  &  2^{-j} \left(\int_{1}^\infty \left(1+2^{j}v^{-2} \right)^{-2N}\gamma(v)^2  v dv \right)^{1/2}  \\
 & \lesssim  \sum_{j\geq 2}  2^{j} \left(\int_{1}^\infty \left(1+2^{j}v^{-2} \right)^{-2N}\gamma(v)^2 v dv \right)^{1/2}.
\end{align*}
To estimate the integral, we decompose for $v\in[1,2^{j/2}]$ and for $v\in[2^{j/2},\infty)$ and we use $\gamma(v)\lesssim (1+v)^{-4}$ to have that
\begin{align*}
\sum_{j\geq 2} \frac{\mu(2^{j+k+1}Q)}{\mu(2^{k+1}Q)} \frac{\nu(2^{k}B)}{\nu(2^{k+j+1}B)} <\infty.
\end{align*}
For the third term of $\alpha_{j,k}$, we have
 \begin{align*}
\sum_{j\geq 2} \frac{\mu(2^{j+k+1}Q)}{\mu(2^{k+1}Q)} \frac{\mu(2^{k}Q)}{\mu( 2^{j+k}Q)2^{2j}} \lesssim \sum_{j\geq 2} 2^{-2j} <\infty.
\end{align*}
We have the desired property due to the additionnal factor $2^{-2j}$, which is obtained by the time-regularity of the semigroup in the case $2-)$. So $(\ref{hypap})$ is satisfied. \findem

\section{Other results.}

\subsection{Maximal regularity on $L^p$ for $p\geq 2$.}

\mb We have the same result for the adjoint operator $T^*$~:
\begin{thm} \label{theo2ap2}  Let $L$ be a generator of a bounded analytic semigroup on $L^2(Y)$. Assume that $(e^{tL^*})_{t>0},(tL^*e^{tL^*})_{t>0}$ and $(t^2L^{2*}e^{tL^*})_{t>0}$ belong to \linebreak[4] ${\mathcal O}_{4}(L^2-L^2)$. Then $T^{*}$ is $H^1_{F,\epsilon,mol}-L^1$ bounded for every $\epsilon>0$, with the Hardy space $H^1_{\epsilon, mol}:=H^1_{\epsilon,mol,{({B_Q}^*)}_{Q\in {\mathcal Q}}}$ (which is the Hardy space constructed with the dual operators $B_Q^*$).
\end{thm}

\deme The adjoint operator $T^*$ is given by~:
$$ T^* f(t,x) = \int _{s=t}^{Z}  \left[L^{*} \left(e^{(s-t) L}\right)^{*} f(s,.)\right](x) ds.$$
The parameter $Z$ depends on the time interval $J$, it is defined by~:
$$ Z := \left\{ \begin{array} {ll}
    \infty & \textrm{ if $J=(0,\infty)$} \\
    l & \textrm{ if $J=(0,l)$}
                \end{array} \right. .$$
The argument of the previous theorem can be repeated and we omit details.
\findem

\mb So now we can apply our general theory to obtain the following result~:

\begin{thm} Let $L$ be a generator of a bounded analytic semigroup on $L^2(Y)$ such that $(e^{tL^*})_{t>0},(tL^*e^{tL^*})_{t>0}$ and $(t^2L^{2*}e^{tL^*})_{t>0}$ belong to ${\mathcal O}_{4}(L^2-L^2)$. Let us assume that for $q_0\in(2,\infty]$, for all balls $Q\subset X$ and all functions $f\in L^2(X)$, we have
$$\left(\frac{1}{\nu(Q)}\int_{Q} \left|A_Q(f)\right|^{q_0} d\mu \right)^{1/q_0} \lesssim \inf_{Q} M_{HL,2}(f).$$
Then for all $p\in (q_0',2]$ the operator $T^*$ is $L^p(X)$-bounded and so $T$ is $L^{p'}$-bounded. We have also the maximal regularity on $L^p(Y)$ for all $p\in[2,q_0)$.
\end{thm}

\deme We use Theorem \ref{theogeneh} for the operator $T^*$
whose hypotheses are satisfied thanks to Theorem \ref{theo2ap2}. \findem

\subsection{Study of our Hardy spaces.}

\mb To finish this paper, we show some results on our Hardy space. First we
have the off-diagonal decay (\ref{decay}).

\begin{prop} Assume that $(e^{tL})_{t>0} \in {\mathcal O}_{p}(L^2-L^2)$ for an integer $p$. For $B_Q$ defined by (\ref{bq}) and (\ref{opA}), we have that for all balls $Q \subset X$
\be{decay2} \forall i\geq 0,\ \forall k\geq 0,\ \forall f\in L^2(2^{k}Q), \qquad \left\| B_Q(f)
\right\|_{2,S_i(2^{k}Q)} \leq C 2^{-M''i} \|f\|_{2,2^{k}Q} \ee with the exponent $M''=\delta/2-1+p$.
\end{prop}

\deme By definition we have just to prove the decay for the $A_Q$ operator. Let $r$ be the radius of $Q$. As previously, we write $Q=I \times B$ where $I$ is an interval of length $r^2$ and $B$ is a ball in $Y$ of radius $r$. Recall that
$$  A_Q(f)(t,x):=\int_{\sigma=0}^{+\infty} \left[\varphi_{r^2} (t-\sigma)  e^{r^2 L} f(\sigma,.) \right](x) d\sigma.$$
For $i\leq 1$, we just use the $L^2(Y)$-boundedness of $A_Q$ to prove (\ref{decay2}). Then for $i\geq 2$ and $(\sigma,y)\in 2^kQ$ if $(t,x)\in
S_i(2^kQ)$ we have that $d((x,t),(\sigma,y))\simeq 2^{k+i}r$ and by using the definition of the parabolic quasi-distance, we
conclude that either $x\in S_i(2^kB)$ either $t\in S_{2i} (2^{2k}I)$. We will study the two
cases~: \\
First for $x\in S_i(2^k B)$, by the off-diagonal estimate
(\ref{offdiagonallow}) we have the estimate : for all $\sigma>0$
\begin{align*}
 \left\| e^{r^2L}(f(\sigma,.)) \right\|_{2,S_i(2^k B)} \lesssim \frac{\nu(2^kB)}{\nu(2^{i+k}B)} \gamma \left(2^{i+k}\right) \left(\frac{\nu(2^{i+k}B)}{\nu(2^kB)} \right)^{1/2} \|f(\sigma,.)\|_{2,2^{k}B}.
\end{align*}
So by the Minkowski inequality, we obtain
\begin{align*}
\lefteqn{\left\| A_Q(f)(t,.) \right\|_{2,S_i(2^kB)} } & & \\
 & & \lesssim \int_{\sigma=0}^{+\infty} \left(1+ \frac{|t-\sigma|}{r^2} \right)^{-N} \frac{\nu(2^kB)}{\nu(2^{i+k}B)} \gamma\left(2^{i+k}\right) \left(\frac{\nu(2^{i+k}B)}{\nu(2^kB)} \right)^{1/2} \|f(\sigma,.)\|_{2,2^{k}B} \frac{d\sigma}{r^2} \\
  & & \lesssim  \left(\frac{\nu(2^kB)}{\nu(2^{i+k}B)}\right)^{1/2} \gamma\left(2^{i+k}\right) \|f\|_{2,2^{k}Q} \frac{1}{r}.
\end{align*}
Then we integrate for $t\in 2^{2(i+k)}I$ to have
$$ \left\| A_Q(f)\right\|_{2,2^{2(i+k)}I \times S_i(2^kB)} \lesssim \left(\frac{\nu(2^kB)}{\nu(2^{i+k}B)}\right)^{1/2} 2^{i+k} \gamma\left(2^{i+k}\right) \|f\|_{2,2^{k}Q}.$$
For the second case, we have $|t-\sigma|\simeq 2^{2(i+k)}r^2$. By using the $L^2(Y)$-boundedness of the semigroup
$$ \left\| e^{r^2L}(f(\sigma,.)) \right\|_{2,2^{i+k}B} \lesssim  \|f(\sigma,.)\|_{2,2^{k}B}.$$
So by the Minkowski inequality, we obtain
\begin{align*}
  \left\| A_Q(f)(t,.) \right\|_{2,2^{k+i}B} & \lesssim \int_{\sigma \in 2^{k} I} \left(1+ 2^{2(k+i)} \right)^{-N}  \|f(\sigma,.)\|_{2,2^{k}B} \frac{d\sigma}{r^2} \\
   & \lesssim  2^{-2(k+i)(N-1)} \|f\|_{2,2^{k}Q} \frac{1}{r}.
\end{align*}
So we can conclude that
$$ \left\| A_Q(f)\right\|_{2,S_{2i}(2^{2k}I) \times 2^{i+k}B} \lesssim 2^{-(N-2)(k+i)} \|f\|_{2,2^{k}Q}.$$
With these two cases, we can conclude (for $N$ any large enough integer)
$$ \left\| A_Q(f)\right\|_{2,S_i(2^k Q)} \lesssim \left(2^{-(N-2)i} + \left(\frac{\nu(2^kB)}{\nu(2^{i+k}B)}\right)^{1/2} 2^{i+k} \gamma\left(2^{i+k}\right) \right) \|f\|_{2,2^{k}Q}$$
which with the decay of $\gamma$ permits to prove the result. \findem

\mb With this decay $M''>\frac{\delta+2}{2}$ (if $p=4$), we have shown that the Hardy spaces
$H^1_{ato}(X)$ and $H^1_{\epsilon, mol}(X)$ are included into the space
$L^1(X)$ (see Proposition \ref{contL1}). \\
In fact we can improve this result, by comparing it with
the classical Hardy space of Coifman-Weiss on $X$.

\begin{prop} \label{aa} Let $\epsilon>0$. The inclusion $H^1_{ato}(X) \subset H^1_{\epsilon, mol}(X) \subset H^1_{CW}(X)$ is equivalent to the fact for all $r>0$, $(e^{rA})^* ({\bf 1}_Y)={\bf 1}_{Y}$ (in the sense of Proposition \ref{inclus}).
\end{prop}

\deme We use the notations of Proposition \ref{inclus}. By using this Proposition, we know that $H^1_{\epsilon, mol}(X) \subset H^1_{CW}(X)$ is equivalent to the fact that for all balls $Q$ of $X$, $A^{*}({\bf 1}_{X}) ={\bf 1}_{X}$ in the sense of $(Mol_{\epsilon,Q})^{*}$. Let $Q=B((t_Q,c_Q),r_Q)$ be fixed. By (\ref{Aetoile}) we know that
$$ A_Q^*(g)(\sigma,x):= \int_{t\in \R^+}  \varphi_{r_Q^2}(t-\sigma) \left[\left(e^{r_Q^2L}\right)^*g(t,.)\right](x) dt.$$
As $\int_\R  \varphi(t) dt =1$, we formally obtain
$$ A_Q^{*}({\bf 1}_X)(\sigma,x) = (e^{r_Q^2L})^{*}({\bf 1}_Y)(x).$$
This equality can be rigorously verified by defining $(e^{r_Q^2L})^{*}({\bf 1}_Y)(x)$ as the continuous linear form on the space
$$ Mol_{\epsilon,r_Q}(Y):= \left\{ f\in L^1(Y),\ \|f\|_{Mol_{\epsilon,r_Q}(Y)}<\infty \right\},$$
where
$$ \|f\|_{Mol_{\epsilon,r_Q}(Y)}:=\sup_{i\geq 0} \|f\|_{2,S_i(Q_Y)} \left(\nu(2^{i}Q_Y)\right)^{1/2} 2^{\epsilon i}.$$
Here we write
$$Q_Y= B(c_Q,r_Q) = \left\{ y\in Y, d_Y(x,c_Q)\leq r_Q \right\}$$
the ball in $Y$. Then the equivalence is a consequence of Proposition \ref{inclus}. \findem

\gb In the paper \cite{ABZ}, the authors have shown that with $-L$ equals to the laplacian on $X$ a complete Riemannian manifold with doubling and Poincar\'e inequality, the operator $T$ is bounded on $H^1_{CW}(X)$ (not just bounded into $L^1(X)$). This is a better result than the one here because Proposition \ref{aa} applies (see \cite{ABZ}) so
$$ H^1_{ato}(X) \subset H^1_{\epsilon, mol}(X) \subset H^1_{CW}(X) \subset L^1(X). $$
But the $H^1_{CW}$-boundedness is using stronger hypotheses than ours in a specific situation.

\end{document}